\documentclass{amsart}
\usepackage{amssymb}
\numberwithin{equation}{section}

\newtheorem{prop}{Proposition}
\newtheorem{thm}[prop]{Theorem}

\newtheorem{lem}[prop]{Lemma}
\newtheorem{conj}[prop]{Conjecture}

\theoremstyle{definition}
\newtheorem{ex}[prop]{Example}
\newtheorem{rem}[prop]{Remark}

\begin{document}

\title[]{Multi-atoms and monotonicity of
generalized Kostka polynomials}

\author{Mark Shimozono}
\address{Dept.\ of Mathematics \\ Virginia Tech \\
Blacksburg, VA}
\email{mshimo@math.vt.edu}

\begin{abstract} The Poincar\'e polynomials of isotypic components
of the graded $GL(n)$-modules given by twists by line bundles
of coordinate rings of closures of conjugacy classes of nilpotent matrices,
are $q$-analogues of the Littlewood-Richardson coefficients
giving multiplicities in a tensor product of irreducible
$GL(n)$-modules indexed by rectangular partitions.  These polynomials
are the $q$-enumeration of a set of Young tableaux (called LR tableaux)
with a generalized charge statistic.  These polynomials satisfy a
monotonicity property which extends that of the Kostka-Foulkes.  In this
paper the monotonicity property is realized combinatorially by a
directed system of embeddings of graded posets whose objects are
sets of LR tableaux.  The final object in this directed system
is a weakened dual of the cyclage poset of Lascoux and Sch\"utzenberger
on the set of column strict tableaux of a fixed content.
For certain sequences of rectangles, the image of the embedding of a set
of LR tableaux into column-strict tableaux, is shown to be a set of
catabolizable tableaux, proving a conjecture of the author and J. Weyman.
\end{abstract}

\maketitle

\newcommand{\key}{\mathrm{key}}
\newcommand{\atom}{\mathrm{matom}}
\newcommand{\content}{\mathrm{content}}
\newcommand{\charge}{\mathrm{charge}}
\newcommand{\cocharge}{\mathrm{cocharge}}

\newcommand{\Z}{\mathbb{Z}}
\newcommand{\la}{\lambda}
\newcommand{\wh}{\widehat{w}}\
\newcommand{\Rhat}{\widehat{R}}
\newcommand{\K}{K}
\newcommand{\tK}{\widetilde{K}}
\newcommand{\Ph}{\widehat{P}}
\newcommand{\Qh}{\widehat{Q}}
\newcommand{\Th}{\widehat{T}}
\newcommand{\Vh}{\widehat{V}}
\newcommand{\NE}{\mathrm{NE}}
\newcommand{\SW}{\mathrm{SW}}
\newcommand{\RRR}{\mathcal{R}}
\newcommand{\NN}{\mathbb{N}}
\newcommand{\ZZ}{\mathbb{Z}}
\newcommand{\C}{\mathbb{C}}
\newcommand{\LRC}{c}
\newcommand{\LRT}{\mathrm{LRT}}
\newcommand{\CST}{\mathrm{CST}}
\newcommand{\rows}{\mathrm{rows}}
\newcommand{\ST}{\mathrm{ST}}
\newcommand{\word}{\mathrm{word}}
\newcommand{\cword}{\mathrm{cword}}
\newcommand{\cmtype}{\mathrm{ctype}}
\newcommand{\NC}{\mathcal{N}}
\newcommand{\SSS}{\mathcal{S}}
\newcommand{\TTT}{\mathcal{T}}
\newcommand{\LLT}{c}
\newcommand{\sign}{\mathrm{sign}}
\newcommand{\weight}{\mathrm{weight}}
\newcommand{\Roots}{\mathrm{Roots}}
\newcommand{\inner}[2]{\langle #1\,,\,#2\rangle}
\newcommand{\sq}{\times}
\newcommand{\Knuth}{\sim_K}
\newcommand{\pr}{\mathrm{pr}}
\newcommand{\rev}{\mathrm{rev}}
\newcommand{\ev}{\mathrm{ev}}
\newcommand{\std}{\mathrm{std}}
\newcommand{\cstd}{\mathrm{cstd}}
\newcommand{\tr}{\mathrm{tr}}
\newcommand{\cod}{{\widetilde{d}}}
\newcommand{\dom}{\trianglerighteq}
\newcommand{\etahat}{\widehat{\eta}}
\newcommand{\hattheta}{\widehat{\theta}}
\newcommand{\colcat}{\mathrm{colcat}}
\newcommand{\cat}{\mathrm{cat}}
\newcommand{\ccat}{\mathrm{cat}}
\newcommand{\Sh}{\widehat{S}}
\newcommand{\Img}{\mathrm{Im}\,\,}

\section{Introduction}

Let $X_\mu\subset gl(n,\C)$ be the Zariski closure
of the conjugacy class of the nilpotent $n\times n$ Jordan matrix with
block sizes given by the parts of the transpose $\mu^t$ of the partition
$\mu$ of $n$.  $X_\mu$ is a subvariety of the affine space $gl(n,\C)$,
so there is a canonical surjection of coordinate rings
$\C[x_{ij}:1\le i,j\le n]=C[gl(n,\C)]\rightarrow \C[X_\mu]$.
Moreover, $X_\mu$ is closed under multiplication by scalar matrices,
so $\C[X_\mu]$ is not only filtered but graded.
The variety $X_\mu$, being the closure of a conjugacy class,
is stable under conjugation, inducing a rational graded action of
$GL(n,\C)$ on $\C[X_\mu]$.

These adjoint orbit closures satisfy $X_\mu\subset X_\nu$ if $\mu\dom\nu$,
inducing surjections $\C[X_\nu]\rightarrow \C[X_\mu]$ of graded
$GL(n)$-modules, given by restriction of functions.  For the
dominant integral weight $\tau$ of $GL(n)$, let
$\K_{\tau;\mu}(q)$ denote the Poincar\'e polynomial of the $\tau$-th
isotypic component of $\C[X_\mu]$.  It follows immediately that the
Poincar\'e polynomials satisfy the monotonicity property
\begin{equation} \label{coord mono}
\K_{\tau;\mu}(q) \le \K_{\tau;\nu}(q)\qquad\text{ if $\mu\dom\nu$}
\end{equation}
(coefficientwise), for every $\tau$.
The (normalized) Kostka-Foulkes polynomial $\tK_{\la,\mu}(q)$ occurs as
the Poincar\'e polynomial of an isotypic component of $\C[X_\mu]$ \cite{W}.
More precisely, let $\{ \omega_i: 1\le i\le n-1 \}$ be the fundamental
weights for $sl_n$.  If $\tau=\sum_{i=1}^{n-1} (\la_i-\la_{i+1}) \omega_i$
for a partition $\la$ of $n$, then
$\K_{\tau;\mu}(q)=\tK_{\la^t,\mu}(q)$ where $\la^t$ is the transpose
of the partition $\la$ \cite{W}.  In this case \eqref{coord mono}
specializes (after transposing $\la$) to
\begin{equation} \label{Kostka mono}
\tK_{\la;\mu}(q) \le \tK_{\la;\nu}(q)\qquad\text{ if $\mu\dom\nu$}
\end{equation}
for every partition $\la$ of $n$.

In \cite{La} \cite{LS1} \cite{LS2} \cite{LS3} Lascoux and
Sch\"utzenberger provided a beautiful combinatorial explanation
for \eqref{Kostka mono}, giving the following results and constructions.
\begin{enumerate}
\item $\tK_{\la,\mu}(q)$ is the $q$-enumeration of the set
$\CST(\la,\mu)$ of column-strict tableaux of shape $\la$ and content
$\mu$ by the map $\cocharge:\CST(\la,\mu)\rightarrow\NN$.
\item The set $\CST(\mu)=\bigcup_\la \CST(\la,\mu)$ has the structure
of a poset (called the cyclage poset) that is graded by cocharge.
\item If $\mu\dom\nu$ then there is an embedding of graded posets
\begin{equation*}
  \theta^\nu_\mu:\CST(\mu)\rightarrow \CST(\nu)
\end{equation*}
that preserves the shapes of tableaux.
\item A map from standard tableaux to partitions is defined (called
catabolism type in \cite{SW}), and it is asserted that the image
of the embedding $\theta^{(1^n)}_\mu$ is the set of standard tableaux
whose catabolism type dominates $\mu$.
\end{enumerate}

The coordinate rings $\C[X_\mu]$ are examples of a larger family
of $\C[X_\mu]$-modules that afford a compatible
graded $GL(n)$-module structure.  Consider the category of
finitely generated graded $\C[x_{ij}]$-modules supported in $X_\mu$
that afford a compatible graded action of $GL(n)$.  In \cite{KKSW}
a set of generators is given for
the graded Grothendieck group of this category.
Among these modules are a distinguished subfamily of modules $M_R$
that may be indexed
by sequences $R$ of rectangular partitions; these are twists of the
coordinate ring $\C[X_\mu]$ by line bundles.  The Poincar\'e
polynomial $\K_{\la;R}(q)$ of the $\la$-th isotypic components of the
module $M_R$ is a $q$-analogue of the Littlewood-Richardson (LR)
coefficient given by the multiplicity of an irreducible $GL(n)$-module
in the tensor product of irreducibles indexed by the rectangular partitions
of $R$.  In \cite{SW} it is observed that there are natural surjections
between certain of these modules, yielding monotonicity
results that generalize \eqref{coord mono}.

In \cite{S} generalizations of 1 and 2 were proven for the Poincar\'e
polynomials $\K_{\la;R}(q)$.  
In this paper a generalization of 3 conjectured in \cite{KS}
is proven using the generalized cyclage poset of
\cite{S}, yielding a combinatorial proof of the monotonicity
property for $\K_{\la;R}(q)$.  In the case that the sequence
of rectangles $R$ is linearly ordered by containment, the
generalization of 4 conjectured in \cite{SW} is proven, obtaining
a new formula for $\K_{\la;R}(q)$ as the $q$-enumeration of
$R$-catabolizable tableaux of shape $\la$ by the charge map of
Lascoux and Sch\"utzenberger.  

A proof is given for the conjecture \cite{KS} that relates the
polynomials $\K_{\la;R}(q)$ and $\K_{\la^t;R^t}(q)$, where $\la^t$ is
the transpose of the partition $\la$ and $R^t$ is the sequence of
transposes of the rectangles in $R$.  Using this it is shown precisely
how our results generalize 3 and 4 above.

All undefined notation and definitions are taken from \cite{S}.

While finishing this paper, the author was informed of the work
of Schilling and Warnaar \cite{ScW}, which has considerable
overlap with this paper and \cite{S} \cite{S2}.

\section{Monotonicity via embeddings}

\subsection{The Poincar\'e polynomials}

Let $R=(R_1,R_2,\dots,R_t)$ be a sequence of rectangular partitions
and $\la$ a partition.  There is a polynomial $\K_{\la;R}(q)$
with integer coefficients such that
$\K_{\la;R}(1)$ is the multiplicity of the $GL(n)$-irreducible
$V_\la$ of highest weight $\la$ in the tensor product
$V_{R_1} \otimes V_{R_2} \otimes \dots \otimes V_{R_t}$ \cite{SW}.
There is a set $\LRT(\la;R)$ of so-called Littlewood-Richardson
(LR) tableaux and a map $\charge_R:\LRT(\la;R)\rightarrow\NN$
(see \cite{KS} and \cite[Sections 2.4,2.6]{S}) such that:

\begin{thm} \cite[Theorem 11]{S} Let $R$ be a dominant sequence of rectangles,
that is, one in which the sequence of numbers of columns of the rectangles
weakly decreases.  Then
\begin{equation*}
  \K_{\la;R}(q) = \sum_{T\in\LRT(\la;R)} q^{\charge_R(T)}.
\end{equation*}
\end{thm}

\subsection{The inequality}
Let $R$ be a sequence of rectangles.  For $k\ge 1$ let
$\xi^k(R)$ be the partition whose parts are the heights of the
rectangles in $R$ of width $k$.  Write $\xi(R)=(\xi^1(R),\xi^2(R),\dots)$.
Observe that $R$ and $R'$ are reorderings of each
other if and only if $\xi(R)=\xi(R')$.
If $\xi=(\xi^1,\xi^2,\dots)$ and $\phi=(\phi^1,\phi^2,\dots)$
are two sequences of partitions, write $\xi\dom\phi$ for the
partial order given by the direct product of dominance partial orders,
that is, if $\xi^i\dom\phi^i$ for all $i$.

Write $R\dom R'$ if $\xi(R)\dom\xi(R')$.  This is a pseudo order;
it is the direct product of dominance partial orders on the
sequences of rectangles modulo reordering.

\begin{thm} \label{rect mono} \cite{SW}
Let $R$ and $R'$ be dominant sequences of rectangles
such that $R\dom R'$.  Then
\begin{equation*}
  \K_{\la;R}(q) \le \K_{\la;R'}(q)
\end{equation*}
coefficientwise.
\end{thm}

This was proved by algebro-geometric methods.  There is a graded
$GL(n)$-module $M_R$ whose $\la$-th isotypic component has
Poincar\'e polynomial given by $\K_{\la;R}(q)$.  Theorem
\ref{rect mono} follows immediately from the
existence of canonical surjections of graded $GL(n)$-modules
$M_{R'}\rightarrow M_R$.

In particular if $R$ and $R'$ are dominant and $\xi(R)=\xi(R')$ then
$\K_{\la;R}(q) = \K_{\la;R'}(q)$.

\subsection{Embeddings}
It is shown in \cite{S} that the set $\LRT(R)=\cup_\la \LRT(\la;R)$
has the structure of a graded poset with order called the $R$-cocyclage
(defined in \cite[Section 3.3]{S}) and grading function given by a
generalized charge map called $\charge_R$ (defined in \cite[Section 2.6]{S}).
We wish to realize Theorem \ref{rect mono}
by exhibiting grade- and order-preserving embeddings
$\theta^{R'}_R:\LRT(\la;R)\rightarrow\LRT(\la;R')$.
Let us recall the constructions of the embeddings defined
in \cite{KS}, which were conjectured to have such properties.

The pseudo order $R \dom R'$ is generated by relations of two forms:
\begin{enumerate}
\item[(E1)] $R\dom R'$ where $R_i=R'_i$ for $i>2$,
$R_1=(k^a)$, $R_2=(k^b)$, $R_1'=(k^{a-1})$, $R_2'=(k^{b+1})$ for some
$a-1\ge b+1$ and $k$.
\item[(E2)] $R\dom \tau_p R$ where $\tau_p R$ is obtained from $R$
by exchanging the $p$-th and $(p+1)$-st rectangles.
\end{enumerate}

For the relation (E2), it is known that the
bijection $\tau_p: \LRT(\la;R)\rightarrow \LRT(\la;R')$
\cite[Section 2.5]{S} is an isomorphism of posets that preserves the
grading \cite[Theorem 20]{S}.  Write
$\theta^{\tau_p R}_R:\LRT(\la;R)\rightarrow\LRT(\la;\tau_p R)$
for this isomorphism.

For relations of the form (E1), let us define an embedding
$\theta^{R'}_R: \LRT(\la;R)\rightarrow\LRT(\la;R')$.
Suppose first that $R=(R_1,R_2)$.  There is a unique embedding
$\iota_{k,\eta_1,\eta_2}:\LRT(\la;R)\rightarrow \LRT(\la;R')$ since
$0\le |\LRT(\la;R)|\le |\LRT(\la;R')| \le 1$; these tableaux 
are described explicitly in \cite[Prop. 33]{S}.  Moreover
in this two-rectangle case
\begin{equation*}
  \charge_{R'}(\iota_{k,\eta_1,\eta_2}(T)) = \charge_R(T)
\end{equation*}
since $\iota_{k,\eta_1,\eta_2}$ preserves shape and both $R$ and $R'$
have the same maximum column size $k$.

\begin{ex} Let $n=6$, $k=3$, $\eta_1=4$, $\eta_2=2$, and $\la=(543321)$.
The tableaux $T\in \LRT(\la;R)$ and $T'\in\LRT(\la;R')$ are given below.
\begin{equation*}
  T=\begin{matrix}
    1&1&1&5&5\\
    2&2&2&6& \\
    3&3&3& & \\
    4&4&4& & \\
    5&6& & & \\
    6& & & & 
  \end{matrix} \qquad
  T'=\begin{matrix}
    1&1&1&4&4\\
    2&2&2&5& \\
    3&3&3& & \\
    4&5&6& & \\
    5&6& & & \\
    6& & & & 
  \end{matrix}
\end{equation*}
\end{ex}

In general, let $R=((k^{\eta_1}),(k^{\eta_2}),\Rhat)$,
$R'=((k^{\eta_1-1}),(k^{\eta_2+1}),\Rhat)$ and $B=[\eta_1+\eta_2]$.
Define $\theta^{R'}_R: \LRT(R)\rightarrow\LRT(R')$ by:
\begin{equation*}
\begin{split}
  \theta^{R'}_R(T)|_B &= \iota_{k,\eta_1,\eta_2}(T|_B) \\
  \theta^{R'}_R(T)|_{[n]-B} &= T|_{[n]-B}
\end{split}
\end{equation*}
where $T|_B$ denotes the restriction of the tableau $T$ to the
subalphabet $B$.  To extend these functions to words, recall
the definition of the set $W(R)$ of $R$-LR words \cite[Section 2.3]{S}.  Let
\begin{equation} \label{iota def}
\begin{split}
\theta^{R'}_R: W(R)&\rightarrow W(R') \\
P(\theta^{R'}_R(w))&=\theta^{R'}_R(P(w)) \\
Q(\theta^{R'}_R(w))&=Q(w) 
\end{split}
\end{equation}
where $P$ is the Schensted $P$ tableau and
$Q$ is the standard row insertion recording tableau.

The maps $\theta^{R'}_R$ and $\tau_p$ shall be called elementary embeddings.

Now let $R\dom R'$ be arbitrary.
Let $R=R^0\dom R^1\dom \dots \dom R^m=R'$ such that for all $i$,
$R^i\dom R^{i+1}$ is a relation of the form (E1) or (E2).  Define an
embedding $\theta^{R'}_R:\LRT(\la;R)\rightarrow\LRT(\la;R')$ by
the composition of elementary embeddings
\begin{equation*}
  \theta^{R'}_R = \theta^{R^m}_{R^{m-1}} \circ \dots
	\circ \theta^{R^2}_{R^1} \circ \theta^{R^1}_{R^0}.
\end{equation*}

\begin{thm} \label{embedding thm}
\begin{enumerate}
\item Let $R\dom R'$.  Then
the map $\theta^{R'}_R$ is independent of the chain in the
pseudo-order from $R$ to $R'$.
\item If $R\dom R'\dom R''$ then
$\theta^{R''}_R = \theta^{R''}_{R'}\circ \theta^{R'}_R$.
\item For $R\dom R'$, $\theta_R^{R'}$ is an injective map
$\LRT(R)\rightarrow \LRT(R')$ that preserves the shape of tableaux
and gives an isomorphism of the $R$-cocyclage poset with
the full subposet induced by its image under $\theta_R^{R'}$.
\item For $R\dom R'$,
\begin{equation*}
  \charge_{R'} \circ \theta_R^{R'} = \charge_R.
\end{equation*}
\end{enumerate}
\end{thm}

Theorem \ref{embedding thm} is proven over the next several sections.

\subsection{Embeddings and $R$-cocyclage}

It is shown that the elementary embedding associated to a
pseudo-order covering relation of the form (E1),
preserves the partial orders on LR tableaux.
The notation is taken from \cite[Section 2.2]{S}.

\begin{prop} \label{embed cover} Let $R\dom R'$ be as in (E1).
Then the elementary embedding
$\theta=\iota_{k,\eta_1,\eta_2}:\LRT(R)\rightarrow\LRT(R')$
is an isomorphism of $\LRT(R)$ under
$R$-cocyclage, with the full subposet given by its image inside
$\LRT(R')$ under $R'$-cocyclage.
\end{prop}
\begin{proof} It must be shown that for all $S,T\in\LRT(R)$,
$S<_R T$ is a covering relation in the $R$-cocyclage poset if and only
if $\theta(S) <_{R'} \theta(T)$ is a covering relation in the
$R'$-cocyclage poset.  Write $S'=\theta(S)$ and $T'=\theta(T)$.

First it is shown that the forward direction suffices.
Let $S'<_{R'} T'$ be an $R'$-cocyclage covering relation with starting
cell $s$.  Since $T$ has the same shape as $T'$, it admits an $R$-cocyclage
covering relation $S" <_R T$ that starts at the same cell $s$.
By the forward direction,
$\theta(S") <_{R'} T'$ is a covering relation that starts at $s$.  This
implies that $\theta(S")=S'$ since both are $R'$-cocyclages of $T'$ starting
at the cell $s$.  But $S'=\theta(S)$ and $\theta$ is injective, so $S"=S$
and the covering relation $S' <_{R'} T'$ is the image under $\theta$ of a
covering relation $S <_R T$.

Let $ux\in W(R)$ such that $T=P(ux)$,
$S=P(\chi_R(ux))=P((w_0^R x)(w_0^R u))$ (where $\chi_R$ and $w_0^R$ are
defined in \cite[Section 3.2]{S}), and the
difference of the shapes of $T$ and $U=P(u)$ is a cell $s$
that lies strictly east of the $a$-th column, where
$a$ is the maximum number of columns among the rectangles in $R$.
Let $A_i$ and $A_i'$ be the subintervals corresponding to the
rectangles $R_i$ and $R_i'$ \cite[Section 2.2]{S}.
By \cite[Theorem 21 (C2), Prop. 23]{S}, $x\in A_i$ for some $i>1$.

It will be shown that $\chi_{R'}(\theta(ux))=\theta(\chi_R(ux))$.
Write $R=(R_1,R_2,\Rhat)$ and $R'=(R_1',R_2',\Rhat)$.
By definition $w_0^R$ acts on $A_j$ by $w_0^{A_j}$ for $1\le j\le 2$
and by $w_0^{\Rhat}$ on $[n]-B$.
Similarly $w_0^{R'}$ acts on $A_j'$ by $w_0^{A_j'}$ for $1\le j\le 2$
and by $w_0^{\Rhat}$ on $[n]-B$.

Write $\theta(ux)=vy$ where $y$ is a letter.

Suppose that $i>2$, in which case $y=x$.
Now $\theta$ affects only the letters of $B$,
and $\chi_R$ only changes the letters outside of $A_i$
by a right shift by one position \cite[Remark 36]{S}.
It follows that $\chi_{R'}(\theta(ux))=\theta(\chi_R(ux))$.

Otherwise $x\in A_2$.  In passing from $ux$ to either
$\chi_{R'}(\theta(ux))$ or $\theta(\chi_R(ux))$,
the letters in $[n]-B$ are merely shifted to the right by one
position.  By restricting to $B$ it may be assumed that $R=(R_1,R_2)$.
By \cite[Prop. 38]{S} the shape of $S=P(\chi_R(ux))$ is obtained
from the shape of $T$ by removing the cell $s$ and adjoining
a cell $s'$ that is uniquely determined by the shape of $T$ and $s$.
Now $P(\theta(ux))=\theta(P(ux))=\theta(T)$ has the same shape as $T$.
The difference of the shapes of $P(u'x')$ and $P(u')$
is the cell of the row insertion recording tableau $Q(u'x')$ that
contains the maximum letter.
But $Q(\theta(ux))=Q(ux)$ so this cell is $s$.
By \cite[Prop. 38]{S} the shape of $P(\chi_{R'}(u'x'))$ is obtained
from that of $P(u'x')$ by removing $s$ and adjoining the same
cell $s'$ as above.  The tableaux $P(\chi_{R'}(\theta(ux)))$
and $P(\theta(\chi_R(ux)))$ are both in $\LRT(R')$ and have
the same shape and $R'$ has only two rectangles,
so the tableaux coincide.  Moreover, since $\theta$ doesn't change
$Q$ tableaux and both $\chi_R$ and $\chi_{R'}$ induce cocyclages
starting at the same corner cell, it follows that
the $Q$ tableaux of the words $\chi_{R'}(\theta(ux))$ 
and $\theta(\chi_R(ux))$ coincide, hence the words themselves do also
by the bijectivity of the RS correspondence.

This shows that $\chi_{R'}(\theta(ux))=\theta(\chi_R(ux))$
Now $P(u')$ has the same shape
as $P(u)$ (since $Q(ux)=Q(u'x')$) and $T$ has the same shape as $T'$.
It follows that $S'=P(\chi_{R'}(\theta(ux)))$ and that
$S' <_{R'} T'$ is a covering relation under $R'$-cocyclage.
\end{proof}

\begin{prop} \label{hat embed}
Suppose $R\dom R'$ such that $R$ and $R'$ are both
dominant.  Say $R_j$ has $\mu_j$ columns for all $j$.
Let $\Rhat$ and $\Rhat'$ be obtained from $R$ and $R'$
by removing all rectangles that have $\mu_1$ columns,
and $B$ the alphabet corresponding to $\Rhat$ (and also $\Rhat'$).
Then $\Rhat\dom\Rhat'$ and for all $R$-cocyclage minimal tableaux
$T\in\LRT(R)$,
\begin{equation} \label{theta hat}
\theta_R^{R'}(T)|_B=\theta_{\Rhat}^{\Rhat'}(T|_B).
\end{equation}
\end{prop}
\begin{proof} Observe that $T\in\LRT(R)$ is $R$-cocyclage minimal
if and only if it has exactly $\mu_1$ columns.  Let $M$ be the sum of
the numbers of rows of rectangles in $R$ (or $R'$) having exactly
$\mu_1$ columns.  Let $B=[M+1,n]$.  Then
by the definition of LR tableau, $T|_{[M]}=Y:=\key((\mu_1^M))$
(see \cite[Section 5.1]{S} for the definition of $\key$).
Let $S=\theta_R^{R'}(T)$.  Since $S$ has the same shape as $T$ it
follows that $S$ is $R'$-cocyclage minimal.   For the same reasons
$S|_{[M]}=Y$.  By definition of LR tableau
$T|_B\in\LRT(\Rhat)$ and $S|_B\in\LRT(\Rhat')$.

By induction it may be assumed that there is a decreasing sequence in the
pseudo order from $R$ to $R'$ that has at most one step of the form (E1).

Suppose first that there is no step of the form (E1).
Without loss of generality it may be assumed that
$\theta_R^{R'}=\tau_p$, which by the dominance assumption must
exchange rectangles with the same number of columns.

If $\tau_p$ exchanges rectangles with $\mu_1$ columns
then it acts by the identity on $R$-cocyclage minimal tableaux
and \eqref{theta hat} holds trivially.  If $\tau_p$ exchanges rectangles
with strictly less than $\mu_1$ columns, then by induction,
\begin{equation*}
  \theta_R^{R'}(T)|_B = \tau_p(T)|_B = \tau_p(T|_B)
	= \theta_{\Rhat}^{\Rhat'}(T|_B)
\end{equation*}
and again \eqref{theta hat} holds.

So it may be assumed that there is a step of the form (E1).  Let $k$ be as
in (E1).  If $k=\mu_1$ then $\theta_R^{R'}$ is the identity on
$R$-cocyclage minimal tableaux and again \eqref{theta hat} holds trivially.

So assume that $k<\mu_1$.  A decreasing chain in the pseudo order
from $R$ to $R'$ can be given by switching the rectangles
$(k^{\eta_1})$ and $(k^{\eta_2})$ to the front, then applying the map
$\iota_{k,\eta_1,\eta_2}$, then switching back to obtain the sequence $R'$.
By \cite[Remark 39]{S}, these rectangle-switching bijections
are easily computed directly and satisfy \eqref{theta hat}.
\end{proof}

\subsection{Embeddings and $\charge_R$}

\begin{prop} \label{embed charge} Let $R\dom R'$.  Then
\begin{equation} \label{theta charge}
  \charge_{R'}(\theta_R^{R'}(T))=\charge_R(T)
\end{equation}
for all $T\in\LRT(R)$.
\end{prop}
\begin{proof} It is enough to assume that $\theta_R^{R'}$ has either
the form $\iota_{k,\eta_1,\eta_2}$ or $\tau_p$ and to show that the left
hand side of \eqref{theta charge} satisfies the axioms (C1) through (C4)
that characterize $\charge_R$ \cite[Theorem 21]{S}.
This is accomplished using Propositions \ref{embed cover} and 
\ref{hat embed}, \cite[Theorem 20]{S}, and induction.
\end{proof}

\subsection{Proof of Theorem \ref{embedding thm}}

\begin{proof} Consider an arbitrary chain
in the pseudo-order from $R$ to $R'$ and consider the composition
$\theta_R^{R'}$ of elementary embeddings defined by this chain.
Then for this map $\theta_R^{R'}$, 
3 holds, in light of Proposition \ref{embed cover} and
\cite[Theorem 20]{S}.  4 also holds by Proposition \ref{embed charge}.

2 follows from 1, for to compute $\theta_R^{R''}$
one may take a chain in the pseudo-order from $R$ to $R'$ and then
from $R'$ to $R''$.

The proof of 1 proceeds by induction on the number of rectangles
in $R$ and on $\charge_R$.  Consider two different chains from $R$
to $R'$ in the pseudo-order and call the corresponding embeddings
$\theta_1$ and $\theta_2$.  Let $T\in\LRT(R)$.  Suppose first
that $T$ is not $R$-cocyclage-minimal.  Then it admits an $R$-cocyclage
covering relation $S<_R T$.  By induction on $\charge(R)$,
$\theta_1(S)=\theta_2(S)$.  But by 3, $\theta_1(T)=\theta_2(T)$.
Otherwise suppose $T$ is $R$-cocyclage-minimal.
By applying some bijections $\tau_p$ to $T$ and then to both
$\theta_1(T)$ and $\theta_2(T)$, it may be assumed that both $R$ and $R'$
are dominant.
Let $a$ be the maximal number of columns among the rectangles
$R_j$.  The proof of Proposition \ref{hat embed} shows that one may
write $T=\Th Y$ and that $\theta_i(T)=\hattheta_i(\Th)Y$ where
$\hattheta_i(\Th)$ computes the map  $\theta_{\Rhat}^{\Rhat'}(\Th)$.
But by induction on the number of rectangles in $R$,
$\hattheta_1(\Th)=\hattheta_2(\Th)$, so $\theta_1(T)=\theta_2(T)$.
\end{proof}

\section{Embedding of LR tableaux into column-strict tableaux}

Given a sequence of rectangles $R$, let $\rows(R)$ be the
sequence of one row rectangles of lengths given by the
weight $\gamma(R)=(\mu_1^{\eta_1},\mu_2^{\eta_2},\dots)\in \ZZ^n$,
that is, $\rows(R)$ is obtained by splitting each of the rectangles
$R_i$ in $R$ into its constituent rows.
Clearly $R\dom\rows(R)$ and $\rows(R)$ is minimal in the
pseudo-order.  Since $\rows(R)$ consists of
one row rectangles, $\LRT(\rows(R))$ is the set of 
column-strict tableaux of content $\gamma(R)$ and partition shape,
with a weak version of the cocyclage poset structure of
Lascoux and Sch\"utzenberger and graded by the usual charge statistic.
But the embedding $\theta_R := \theta^{\rows(R)}_R$ preserves grading
and the shape of tableaux by Theorems \ref{embed cover} and
\ref{embed charge}.  For certain $R$ it is shown that the image of $\theta_R$
is given by the $R$-catabolizable tableaux of \cite{SW}.
Thus $\K_{\la;R}(q)$ has another formula,
as the $q$-enumeration of the $R$-catabolizable tableaux of
shape $\la$ by the charge statistic.

The map $\theta_R$ is a generalization of the
standardization map in \cite{La} which is a
cyclage- and cocharge-preserving embedding 
of the column-strict tableaux into standard tableaux; a precise
statement of this is given in Remark \ref{cyc std}.

\subsection{Catabolizable tableaux}
\label{cat sec}

Given a (possibly skew) column-strict tableau $S$ and index $r$
(resp. index $c$), let $H_r(S) = P(S_n S_s)$
(resp. $V_c(S)=P(S_e S_w)$), where $S_n$ and $S_s$ (resp.
$S_e$ and $S_w$) are the north and south (resp. east and west)
subtableaux obtained by slicing $S$ horizontally
(resp. vertically) between its $r$-th and $(r+1)$-st rows
(resp. $c$-th and $(c+1)$-st columns).

Let $S$ be a column-strict tableau of partition shape in the alphabet $[n]$,
$R=(R_1,R_2,\dots)$ such that $R_j$ is the rectangular partition with
$\eta_j$ rows and $\mu_j$ columns, and $\Rhat=(R_2,R_3,\dots)$.
Suppose $S|_{A_1}=Y_1$.  In this case the
$R_1$-catabolism of $S$ (the $R_1$-column-catabolism of $S$)
is defined to be the tableau $\cat_{R_1}(S)=H_{\eta_1}(S-Y_1)$
(resp. $\colcat_{R_1}(S)=V_{\mu_1}(S-Y_1)$).

Define $R$-catabolizability (resp. $R$-column-catabolizability)
by the following rules.
\begin{enumerate}
\item The empty tableau is the unique (resp. column) catabolizable tableau
for the empty sequence of rectangles.
\item Otherwise, $S$ is $R$-catabolizable (resp. $R$-column-catabolizable)
if and only if $S|_{A_1} = Y_1$ and $\cat_{R_1}(S)$
(resp. $\colcat_{R_1}(S)$) is $\Rhat$-catabolizable
(resp. $\Rhat$-column-catabolizable) in the alphabet $[\eta_1+1,n]$.
\end{enumerate}
Observe that if $S$ is $R$-catabolizable then $\content(S)=\gamma(R)$,
and that if $R$ is dominant and consists of single rows, then
$S$ is $R$-catabolizable if and only if $S$ is column-strict.

\begin{ex} Let $R=((33),(33),(222))$ so that $\mu=(332)$,
$\eta=(223)$, $n=7$, $A_1=[1,2]$, $A_2=[3,4]$, $A_3=[5,7]$, and
\begin{equation*}
  Y_1 = \begin{matrix} 1&1&1\\2&2&2 \end{matrix} \qquad
  Y_2 = \begin{matrix} 3&3&3\\4&4&4 \end{matrix} \qquad
  Y_3 = \begin{matrix} 5&5\\6&6\\7&7 \end{matrix}
\end{equation*}
The following computation shows that the tableau $S$
is both catabolizable and column-catabolizable.
\begin{equation*}
\begin{split}
  S &= \begin{matrix} 1&1&1&3&4&5\\2&2&2&4&5&6\\3&3&6& & & \\4&7&7& & &
  \end{matrix} \\
  S_e S_w &= \begin{matrix}
    \sq&\sq&\sq&3&3&6\\
    \sq&\sq&\sq&4&7&7\\
    3&4&5& & & \\
    4&5&6& & &     
  \end{matrix} \\
  \ccat_{R_1}(S)&=P(S_e S_w) =
  \begin{matrix}
    3&3&3&6&7\\
    4&4&4&7& \\
    5&5& & & \\
    6& & & & 
  \end{matrix} \\
  \cat_{R_2} \cat_{R_1} (S) &= Y_3 \Knuth
  \begin{matrix}
    \sq&\sq&5&5\\
    \sq&\sq&6& \\
    6&7& & \\
    7& & &     
  \end{matrix}
\end{split}
\end{equation*}
\end{ex}

The following result solves a special case of \cite[Conjecture 26]{SW}.

\begin{thm} \label{embed image} Suppose $R$ is dominant
and $R^t$ is dominant, that is, $R_j \supseteq R_{j+1}$ for all $j$.
Then the image of $\theta_R$ is the set of $R$-catabolizable tableaux.
\end{thm}

The rest of the section follows \cite{SW}.
Two weak versions of the cocyclage poset
on column-strict tableaux \cite{LS3} are now defined.
Let $u$ be a word and $x$ a letter,
$S=P(ux)$, $U=P(u)$, $T=P(xu)$, and
$s$ (resp. $s'$) the cell that is the difference of the shapes of $S$
(resp. $T$) and $U$.  Contrary to \cite{SW}, define
$T <_{(,>c)} S$ (resp. $T <_{(>r,)} S$)
if the cell $s=(i,j)$ (resp. $s'=(i',j')$) satisfies $j>c$ (resp. $i'>r$).

Given a sequence of rectangles $R$, let $R_j$ have $\eta_j$ rows
and $\mu_j$ columns, $a=\max_j \mu_j$
and $r>0$.  Then the set $\CST(\gamma(R))$ is
a graded poset under the order $<_{(,>a)}$ and
also under $<_{(>r,)}$.  This follows from
\cite[Lemme 2.13]{LS3} and variants.  Let $b=\max_j \eta_j$.
Only the case $r=b$ is used here.

\begin{prop} \label{col cat cyc}
Let $R$ be dominant, $S$ an $R$-column-catabolizable tableau and
$T <_{(,>a)} S$.  Then $T$ is $R$-column catabolizable.
\end{prop}
\begin{proof} Observe that in the above notation, $S$, $T$, and $U$
all contain $Y_1$.  So define $S_e$ and $S_w$ as in the definition
of $\colcat_{R_1}(S)$ and similarly for
$T$ and $U$.  In particular $S_e = P(U_e x)$ and $S_w=U_w$.
By assumption $S$ is $R$-column-catabolizable,
so $P(S_e S_w)=\colcat_{R_1}(S)$ is $\Rhat$-column-catabolizable where
$\Rhat=(R_2,R_3,\dots)$.  By definition it suffices to show that
$\colcat_{R_1}(T)=P(T_e T_w)$ is $\Rhat$-column-catabolizable.

Suppose that during the column insertion of $x$ into $U$,
that a letter $y$ is bumped from the $a$-th column
of $U$ to the $(a+1)$-st.  This can only happen if $U$ has
$a$ columns.  Then $P(x U_w)=P(T_w y)$ where
the difference of the shapes of $P(T_w y)$ and $T_w$
is a cell in the last column, which is the $(a+1)$-st since
$U_w$ and hence $T_w$ has $a$ columns.  Finally $P(y U_e)=T_e$.  Then
\begin{equation*}
 \colcat_{R_1}(T) = P(T_e T_w) = P(y U_e T_w)
\end{equation*}
and
\begin{equation*}
  P(U_e T_w y)=P(U_e x U_w)=P(S_e S_w)=\colcat_{R_1}(S).
\end{equation*}
Since the row insertion of $y$ into $T_w$ ends at the cell $(a+1,1)$,
by \cite[Lemma 22]{S3} the row insertion of $y$ into
$P(U_e T_w)$ ends strictly to the east of the $a$-th column.  It follows
that
\begin{equation*}
  \colcat_{R_1}(T) = P(U_e T_w y) <_{(,>a)} P(y U_e T_w) =
  \colcat_{R_1}(S).
\end{equation*}  
By induction $\colcat_{R_1}(T)$ is $\Rhat$-column-catabolizable.

Otherwise suppose that no such $y$ exists; then
$S_e = P(U_e x)$, $S_w=U_w$, $P(x U_w)=T_w$, $T_e=U_e$, and
\begin{equation*}
\begin{split}
  \colcat_{R_1}(T) &= P(T_e T_w) = P(U_e T_w) \\
  &= P(U_e x U_w) = P(S_e S_w)=\colcat_{R_1}(S).
\end{split}
\end{equation*}
\end{proof}

In an analogous manner one proves the following.

\begin{prop} \label{cat cyc}
Suppose $R^t$ is dominant, $T$ is $R$-catabolizable and
$T <_{(>b,)} S$.  Then so is $S$.
\end{prop}

The next result relates the 
catabolism and column-catabolism operators.

\begin{prop} \label{row col cat}
Suppose $R$ is dominant, $R^t$ is dominant,
and $S$ contains $Y_1$.  Then
\begin{enumerate}
\item $\colcat_{R_1}(S) \le_{(,>a)} \cat_{R_1}(S)$.
\item $\colcat_{R_1}(S) \le_{(>b,)} \cat_{R_1}(S)$.
\end{enumerate}
\end{prop}
\begin{proof} Slice $S$ into four parts, cutting both horizontally
as in the definition of $\cat_{R_1}(S)$ and vertically as in
that of $\colcat_{R_1}(S)$.  Denote the northeast, southeast, and
southwest tableaux by $S_{ne}$, $S_{se}$, and $S_{sw}$.
Then $\colcat_{R_1}(S)=P(S_{se} S_{ne} S_{sw})$
and $\cat_{R_1}(S)=P(S_{ne} S_{sw} S_{se})$.

Now the row insertion of $S_{se}$ into $S_{sw}$ creates all cells
strictly east of the $a$-th column.  It follows that the same
happens for the row insertion of $S_{se}$ into $P(S_{ne} S_{sw})$
and that a sequence of $<_{(,>a)}$-cocyclages lead from
$\cat_{R_1}(S)$ to $\colcat_{R_1}(S)$.  This proves 1.
2 is proven similarly.
\end{proof}

For special sequences of rectangles, the two kinds of $R$-catabolizability
are equivalent.

\begin{prop} \label{row col cat equal}
Let $R$ be dominant and $R^t$ dominant.  Then
$S$ is $R$-catabolizable if and only if $S$ is $R$-column-catabolizable.
\end{prop}
\begin{proof} It may be assumed that $S$ contains $Y_1$.
Under this assumption each of the following implies the next,
using induction and Propositions \ref{col cat cyc} and \ref{row col cat}.
\begin{enumerate}
\item $S$ is $R$-catabolizable.
\item $\cat_{R_1}(S)$ is $\Rhat$-catabolizable.
\item $\cat_{R_1}(S)$ is $\Rhat$-column-catabolizable.
\item $\colcat_{R_1}(S)$ is $\Rhat$-column-catabolizable.
\item $S$ is $R$-column-catabolizable.
\end{enumerate}
To reverse the implications, replace 3 by
``$\colcat_{R_1}(S)$ is $\Rhat$-catabolizable" and
Proposition \ref{col cat cyc} by \ref{cat cyc}.  Then
each statement implies the previous one.
\end{proof}

\subsection{Strong $R$-cocyclage}

Our proof of Theorem \ref{embed image} requires
the strongest version of the cocyclage poset.  Suppose $u x\in W(R)$ with
$x$ a letter.  Write $P(\chi_R(ux)) < P(ux)$ if, for every element $\tau$
in the symmetric group that permutes $R$, the last letter of $\tau (ux)$
is not in the first 
subalphabet of $\tau R$; call the resulting partial order the
strong $R$-cocyclage.  By \cite[Prop. 23]{S} this holds if and only if
$\charge_R(ux)=\charge_R(\chi_R(ux))-1$.  The strong $R$-cocyclage
poset is graded by $\charge_R$ and is obtained from the $R$-cocyclage poset
by adding some relations.

\begin{thm} \label{strong embed} Suppose 
some rearrangement of $R$ is decreasing with respect to containment
and $R\dom R'$.  Then $\theta=\theta_R^{R'}$ is
an embedding of the strong $R$-cocyclage poset with
the full subposet of its image in the strong $R'$-cocyclage poset.
Moreover, if $w\in W(R)$ and $P(\chi_R(w))<P(w)$ then
\begin{equation}\label{theta chi}
  \theta(\chi_R(w))=\chi_{R'}(\theta(w))
\end{equation}
\end{thm}

In general it is difficult to check the
strong $R$-cocyclage condition; one must
examine the entire orbit of a word.  However, if the rectangles
of $R$ are linearly ordered by containment, then looking at
one orbit element suffices.

\begin{lem} \label{linear R cocyclage}
Suppose $R_j\supset R_{j+1}$ for all $j$ and $ux\in W(R)$
with $x$ a letter.  Then $P(\chi_R(ux))<P(ux)$ if and only if $x$ is
not in the first subalphabet $A_1$ for $R$.
\end{lem}
\begin{proof} The forward direction holds by definition.
Consider the two rectangle case $R=(R_1,R_2)$.
Write $vy=\tau_1(ux)$ and $A_j$ and $A_j'$ the subalphabets
for $R$ and $R'=(R_2,R_1)$.  Since $R_1\supseteq R_2$, it follows that
$x\in A_1$ if and only if $y\in A_1'$ by \cite[Prop. 38]{S}.

Now let $R=(R_1,\dots,R_t)$.
Since $\tau_p$ for $p>1$ does not change letters in the
first subalphabet, it is enough to check the last letter of
$\tau(ux)$ for $\tau$ a minimal coset representative in
$S_1\times S_{t-1} \backslash S_t$, that is, $\tau$ has the form
$\tau=\tau_1\tau_2\dots \tau_r$.  Now $\tau R$ starts with
$(R_r,R_1,\dots)$.  By restriction to these first two subalphabets and the
two rectangle case above, it follows that 
the last letter of $\tau(ux)$ is not in the first subalphabet.
But this holds for all $r$, so $P(\chi_R(ux))<P(ux)$.
\end{proof}

Proof of Theorem \ref{strong embed}:
\begin{proof} It follows immediately from the definitions that if two
sequences of rectangles are rearrangements of each other, then their
strong cocyclage posets are isomorphic.  It is clear that since
$R$ has a rearrangement that is decreasing with respect to containment,
so must $R'$.  So it may be assumed without loss of generality that $R$ 
and $R'$ are related as in (E1).  Write $A_j$ and $A_j'$ for the subalphabets
for $R_j$ and $R_j'$ respectively.  Let $R=(R_1,\dots,R_t)$.

Let $w=ux$ and $vy=\theta(ux)$ with $x$ and $y$ letters.
By restriction to the first two subalphabets
and the explicit description of the two-rectangle LR tableaux
in \cite[Prop. 33, 38]{S}, it follows that $x\in A_1$ if and only if
$y\in A_1'$.  By Lemma \ref{linear R cocyclage} this establishes Theorem
\ref{strong embed} except for the equation \eqref{theta chi}.

Assume $P(\chi_R(ux))<P(ux)$ so that $x\not\in A_1$ and $y\not\in A_1'$.
Consider the two words in \eqref{theta chi}.
Suppose first that $x\in A_j$ for $j>2$.  By definition, $\theta$
only changes letters in the first two subalphabets, so $y=x$.
By \cite[Remark 39]{S} $\chi_R$ only affect the
alphabets other than $A_j$ by shifting letters to the right by
one position.  In this case it is clear that \eqref{theta chi} holds.
Since $x\not\in A_1$ the only other case is $x\in A_2$ (and hence
$y\in A_2'$).  Again by \cite[Remark 39]{S} the letters in 
$A_j$ for $j>2$ are merely shifted to the right by $\chi_R$, so
by restriction to $A_1\cup A_2$ it may be assumed that $R=(R_1,R_2)$.
In this case \cite[Prop. 38]{S} establishes \eqref{theta chi}.
\end{proof}

\subsection{Proof of Theorem \ref{embed image}}.

\begin{lem} \label{add rectangle}
Let $R$ be a dominant sequence of rectangles, $R=(R_1,\Rhat)$,
and $u\in W(\Rhat)$ in the alphabet $[\eta_1+1,n]$.
Then $\theta_R(u Y_1) = \theta_{\Rhat}(u) Y_1$
where $\theta_{\Rhat}$ is understood to have image in the alphabet
$[\eta_1+1,n]$ and the tableau $Y_1$ is identified with its row-reading word.
\end{lem}
\begin{proof} Observe that $u Y_1\in W(R)$ since $u\in W(\Rhat)$,
so that $\theta_R(u Y_1)$ makes sense.

It suffices to show that the two words
$\theta_R(u Y_1)$ and $\theta_{\Rhat}(u) Y_1$ have the
same $P$ tableaux and same $Q$ tableaux.

Let $M$ be the length of $u$ and $N$ the length of $u Y_1$.
Then by the definition of recording tableau and $\theta$,
\begin{equation*}
\begin{split}
  Q(\theta_{\Rhat}(u) Y_1)|_{[M]} &= Q(\theta_{\Rhat}(u)) \\
  &= Q(u) = Q(u Y_1)|_{[M]} = Q(\theta_R(u Y_1))|_{[M]}.
\end{split}
\end{equation*}
Observe that the letters of $Y_1$ are smaller than those in $u$, so
$P(u Y_1)$ is obtained from $P(u)$ by pushing each of the first $\mu_1$
columns down by $\eta_1$ cells and placing $Y_1$ in the vacated positions.
The same is true for $P(\theta_{\Rhat}(u))$ and
$P(\theta_{\Rhat}(u) Y_1)$.  Moreover it is clear that for all
$1\le j\le \mu_1$, the $j$-th column of $Q(u Y_1)|_{[M+1,N]}$
is equal to the $j$-th column of $Q(Y_1)+M$.  This shows the
equality of recording tableaux
\begin{equation*}
  Q(\theta_R(u Y_1))=Q(u Y_1)=Q(\theta_{\Rhat}(u) Y_1).
\end{equation*}

For equality of $P$-tableaux, let $U=P(u)$ have shape $\la$, say.

Suppose first that $\la_1>\mu_1$.  Let $U' <_{\Rhat} U$ be the
$\Rhat$-cocyclage on $U$ starting at the cell $s$ at the
bottom of the rightmost column $\la_1$.
Since the letters of $Y_1$ are smaller
than those in $U$, it follows that $P(U' Y_1) <_R P(U Y_1)$.  We have
\begin{equation*}
\begin{split}
  P(\theta_R(U Y_1) &= \theta_R(P(U Y_1)) \\
  &>_{\rows(R)} \theta_R(P(U' Y_1)) \\
  &= P(\theta_R(U' Y_1)) \\
  &= P(\theta_{\Rhat}(U') Y_1).
\end{split}
\end{equation*}
by the definition of $\theta_R$, Theorem \ref{embedding thm} and
and induction.

On the other hand,
$\theta_{\Rhat}(U') <_{\rows(\Rhat)} \theta_{\Rhat}(U)$.
by Theorem \ref{embedding thm} and \cite[Prop. 15]{S}.
Applying the above argument to add $\eta_1$ single rows
$(\mu_1)$ to the front of the sequence $\rows(\Rhat)$, it follows that
$P(\theta_{\Rhat}(U')Y_1) <_{\rows(R)} P(\theta_{\Rhat}(U) Y_1)$.

There are two $\rows(R)$-cocyclages to
$P(\theta_{\Rhat}(U')Y_1)$ from either
$P(\theta_R(U Y_1))$ or $P(\theta_{\Rhat}(U) Y_1)$
and both are induced by a cell in the last column
on tableaux of the same shape.  It follows that
$P(\theta_R(U Y_1))=P(\theta_{\Rhat}(U) Y_1)$.

The other case is that $\la_1 \le \mu_1$.
Consider the following way to compute $\theta_R$:
\begin{equation*}
\begin{split}
  (R_1,R_2,\dots) &\dom ((\mu_1)^m,R_2,R_3,\dots) \\
  &\dom (R_2,R_3,\dots,(\mu_1)^m) \\
  &\dom (\rows(\Rhat),(\mu_1)^m) \\
  &\dom \rows(R).
\end{split}
\end{equation*}
Let us describe the corresponding embeddings of sets of LR tableaux.
It is not hard to show that the first embedding, which chops the
first rectangle $R_1=(\mu_1^m)$ into its constituent rows,
is the identity on the appropriate LR tableaux.
The second and fourth maps are compositions of rectangle-switching
bijections.  The third map is equal to $\theta_{\Rhat}$ when restricted
to the alphabet $[n-m]$, and the identity on $[n-m+1,n]$.
Moreover, the two rectangle-switching steps always exchange
subtableaux vertically since the number of columns in the tableau
$P(U Y_1)$ is equal to the number of columns in $Y_1$ \cite[Remark 39]{S}.

This given, tracing the image of the tableau $P(U Y_1)$ through this map
produces the tableaux $P(U Y_1)$, $P((Y_1+(n-m)) (U-m))$,
$P((Y_1+(n-m)) (\theta_{\Rhat}(U)-m))$, and finally
$P(\theta_{\Rhat}(U) Y_1)$.  By the previous computation and the
definition of $\theta$ we have
\begin{equation*}
\begin{split}
  P(\theta_R(u Y_1))&=\theta_R(P(u Y_1))=\theta_R(P(U Y_1)) \\
  &= P(\theta_R(U Y_1)) = P(\theta_{\Rhat}(U) Y_1) \\
  &= P(\theta_{\Rhat}(P(u)) Y_1) = P(P(\theta_{\Rhat}(u)) Y_1) \\
  &= P(\theta_{\Rhat}(u) Y_1).
\end{split}
\end{equation*}
\end{proof}

Proof of Theorem \ref{embed image}:
\begin{proof} The proof proceeds by induction on the number of rectangles
in $R$.  Let $S$ be a column-strict tableau
of shape $\la$ and content $\gamma(R)$.

For the forward direction, suppose $S=\theta_R(T)$ for some $T\in\LRT(R)$.
Define $T'_e = w_0^R T_e$ and $T'_w=w_0^R T_w$.
Then $P(T'_e T'_w Y_1) <_R T$.  One may show that
inducing corresponding cocyclages on $S$ yields
$P(S_e S_w Y_1) <_{(,>a)} S$.  By Theorem \ref{embedding thm},
$\theta_R(T'_e T'_w Y_1) = S_e S_w Y_1$.
By Lemma \ref{add rectangle} this implies
$\theta_{\Rhat}(T'_e T'_w) = S_e S_w$, so that
$\colcat(S)=P(S_e S_w)$ is in the image of $\theta_{\Rhat}$.
By induction $\colcat(S)$ is $\Rhat$-column-catabolizable, so
that $S$ is column-catabolizable by definition.

Conversely suppose $S$ is $R$-column-catabolizable.  Then $S|_{A_1}=Y_1$ and
$\Sh=\colcat_{R_1}(S)=P(S_e S_w)$ is $\Rhat$-column-catabolizable.
By induction there is a tableau $\Th\in\LRT(\Rhat)$ such that
$\theta_{\Rhat}(\Th)=\Sh$.  Since $\theta_{\Rhat}$ preserves recording
tableaux, there are tableaux $T'_e$ and $T'_w$ of the same shapes
as $S_e$ and $S_w$ respectively, such that $P(T'_e T'_w)=\Th$ and 
$\theta_{\Rhat}(T'_e T'_w) = S_e S_w$.  By Lemma \ref{add rectangle}
$\theta_R(T'_e T'_w Y_1)=S_e S_w Y_1$.  Write
$T_e = w_0^R T'_e$ and $T_w = w_0^R T'_w$.  Let $j$ be the number of
letters in $T'_e$.  Observe that one has a sequence of covering
relations in the strong $R$-cocyclage that gives
\begin{equation*}
  P(\chi_R^j(T_w Y_1 T_e)) < P(T_w Y_1 T_e).
\end{equation*}
By Theorem \ref{strong embed},
\begin{equation*}
  \chi_{\rows(R)}^j(\theta_R(T_w Y_1 T_e)) =
  \theta_R(\chi_R^j(T_w Y_1 T_e)) =
  \theta_R(T'_e T'_w Y_1) = S_e S_w Y_1.
\end{equation*}
Now $\chi_{\rows(R)}$ is merely right circular rotation since
$w_0^{\rows(R)}$ is the identity.  Then
\begin{equation*}
  P(\theta_R(T_w Y_1 T_e)) = P(S_w Y_1 S_e) = S
\end{equation*}
which shows that $S$ is in the image of $\theta$.
\end{proof}

\subsection{Multi-atoms and catabolism multi-type}

In this section the notion of an atom \cite{La} is generalized.

Let $\gamma$ be a partition of length $n$.
Consider the collection of (dominant) sequences of rectangles $R$ such that
$\gamma(R)=\gamma$.  For each such $R$ there is an embedding
$\theta_R:\LRT(R)\rightarrow \CST(\gamma)$.
If $R\dom R'$ then $\gamma(R')=\gamma(R)$ and 
the image of $\theta_{R'}$ contains the image of $\theta_R$,
by Theorem \ref{embedding thm}.  Define the multi-atom
$\atom(R)$ to be the set of $S\in\CST(\gamma)$ such that
$S$ is in the image of $\theta_R$ but not in the image of
$\theta_{R'}$ where $R\dom R'$ and $R'\not\dom R$.

We give a somewhat ackward characterization of the
image of $\theta_R$ in $\CST(\gamma(R))$ for arbitrary $R$.
This condition, which requires the calculation of several
images of a given word of content $\gamma$ under the
automorphisms of conjugation, shows that the multi-atoms for $R$ such that
$\gamma=\gamma(R)$, induce a set partition of $\CST(\gamma)$.
The multi-atoms $\atom(R)$ where $R$ is linearly ordered by
containment, have a nicer characterization.  Given a tableau $S$
of content $\gamma$, we construct a sequence of partitions
$\cmtype(S)=(\xi^1(S),\xi^2(S),\dots)$ called its catabolism multi-type.

\begin{conj} \label{atom conj} Let $\gamma$ be a partition
with at most $n$ parts, $S\in \CST(\gamma)$, and $R$ such that
$\gamma(R)=\gamma$.  Then 
$S \in\atom(R)$ if and only if $\cmtype(S)=\xi(R)$.  In particular,
\begin{equation*}
  \K_{\la;R}(q) = \sum_S q^{\charge(S)}
\end{equation*}
where $S$ runs over the column-strict tableaux of shape $\la$,
content $\gamma$, such that $\cmtype(S) \dom \xi(R)$.
\end{conj}

\begin{thm} \label{atom thm} Conjecture \ref{atom conj}
holds for $R$ such that $R_j\supseteq R_{j+1}$ for all $j$.
\end{thm}

\subsection{Image characterization}

Let $\gamma\in\Z^n$ and $R$ be related by
$\gamma=\gamma(R)$.  For a given positive integer $k$,
let $R^k$ be the subsequence of $R$ consisting of the
rectangles that have exactly $k$ columns.

\begin{prop} \label{image}
Let $S\in\CST(\gamma(R))$.  The following are equivalent:
\begin{enumerate}
\item $S \in \Img \theta_R$.
\item For all $k\ge 1$ there is a permutation $\sigma\in S_n$ such
that $\sigma$ sends the union of the subalphabets of the rectangles in
$R^k$, to an initial interval $I\subseteq [n]$ and 
$(\sigma S)|_I\in \Img \theta_{R^k}$ where
$\sigma$ acts by an automorphism of conjugation.
\end{enumerate}
\end{prop}
\begin{proof} Let $S\in\CST(\gamma(R))$.  Suppose that $S$ has strictly
more than $a=\max_j \mu_j$ columns, so that there is a tableau $T$ such that
$T <_{(,>a)} S$.  By Theorem \ref{embedding thm} and induction on
$\charge_R$ it is enough to show that $S$ satisfies 2 if and only if $T$
does.  Write $S=P(Ux)$ with $x$ a letter, $U$ a tableau such that
the difference of the shapes of $S$ and $U$ is a cell $s$ that
lies in a column strictly east of the $a$-th, and $T=P(xU)$.
Fix $k$, $\sigma$, and $I$ as in the hypothesis.
Write $\sigma(Ux)=U'x'$.  Then $P(x'U')<_{(,>a)} P(U'x')$, with the
cocyclage starting at the same cell $s$.  If $x'\in I$ then
since $I$ is an initial interval, the cell $s'$ that comprises the
difference of the shapes of $P(U'x')|_I$ and $U'|_I$, is weakly east
of $s$, so that
\begin{equation*}
\begin{split}
  (\sigma(T))|_I &= \sigma(P(xU))|_I = P(\sigma(xU))|_I \\
  &= P((x'U')|_I) = P(x' (U'|_I)) <_{(,>a)} P((U'|_I) x') \\
  &= P(U'x')|_I = \sigma(S)|_I.
\end{split}
\end{equation*}
This case is completed by applying Theorem \ref{embedding thm}
for $\theta_{R^k}$.

Otherwise $x'\not\in I$.  Then $\sigma(T)|_I=\sigma(S)|_I$,
and the equivalence of 2 for $k$, $T$ and $S$ is obvious in this case.

Otherwise it may be assumed that $S$ has exactly $a$ columns
(and also that $R$ is dominant).
If $k=a$, then $\sigma$ is the identity and $S|_I$, being a
canonical rectangular tableau with $a$ columns, is automatically in the
image of any map $\theta_{R^k}$.
If $k<a$, in light of \cite[Remark 39]{S},
it is not hard to show that $\sigma(S)|_I$
is also obtained by removing the rectangles in $R^a$
by restriction to the complement of an initial interval, followed
by an automorphism of conjugation that moves the union of alphabets
of $R^k$ to an initial interval, then restricts to that initial
interval.  The proof is completed by Lemma \ref{add rectangle} and
induction on the number of rectangles in $R$.
\end{proof}

\subsection{Nilpotent orbit case}

Suppose that $R$ is a sequence of partitions, all of which have
$k$ columns.  Let $S$ have content $(k^n)=\gamma(R)$.
We associate to $S$ a partition $\xi^k$ as follows.
Consider the sequence of tableaux $S$, $V_k(S)$, $V_k^2(S)$, etc.
If $V_k^j(S)$ has more than $k$ columns, then passing from $V_k^j(S)$
to $V_k^{j+1}(S)$ is a composition of $(,>k)$-cocyclages and the
charge drops.  So the sequence stabilizes.  Let
$y_k(S)$ be the positive integer $j$ such that
$S|_{[j]}=\key((k^j))$.  Define $\xi^k_1=y_k(S)$ and
$\xi^k_{j+1}=y_k(V_k^j(S))-y_k(V_k^{j-1})$.

\begin{lem} \label{two cat rect} For $S$ of content $(k^n)$,
$\xi^k(S)$ is a partition.
\end{lem}
\begin{proof} By restriction to the initial subinterval $[1,i]$
where $i$ is maximal such that $\gamma_1=\gamma_i$ and induction,
it may be assumed that $\gamma=(k^n)$.
Let $S$ be column-strict of content $(k^n)$ and partition shape.
Write $y=y_k$, $\xi^k=\xi$ and $m=y(S)$.  Let $Y$ be the
column-strict tableau of shape $(k^m)$ whose $i$-th row
consists of $k$ copies of the letter $i$, for all $1\le i\le m$.
Let $S_e$ and $S_w$ be as in the definition of $\ccat_{(k^m)}(S)$
and write $\Sh=\ccat_{(k^m)}(S)$.  By induction
$\xi(\Sh)$ is a partition of the number $n-m$.
Since $P(\Sh Y)=V_k(S)$, it follows that the $j$-th part of $\xi(\Sh)$
is the $(j-1)$-st part of $\xi(S)$.  This given, it is enough to show
that the second part $m'$ of $\xi(S)$ does not exceed the first part $m$.
For this, without loss of generality, by restriction to the subalphabet
$[m+m']$, it may be assumed that $\xi(S)=(m,m')$.
This means that $P(S_e S_w) = Y'$ where $Y'$ is the tableau of shape
$(k^{m'}$ whose $i$-th column consists of $k$ copies of the letter
$m+i$ for all $1\le i\le m'$.  By \cite[Prop. 32]{S},
the tableau $S_w$ must be lattice in the alphabet $[m+1,m+m']$
and the shape of $S_e$ must be the 180-degree rotation of the
skew shape given by the complement of the shape of $S_w$ in the
rectangular shape $(k^{m'})$.
By the maximality of $m$, $S_w$ cannot contain all $k$ of the
letters $m+1$, so this means that $S_e$ has exactly $m'$ rows.
On the other hand, since $S$ has partition shape,
$S_e$ has at most $m$ rows.  Therefore $m\ge m'$.
\end{proof}

\subsection{Catabolism multi-type}

Let $\gamma$ be a partition.  $\cmtype(S)=(\xi^1(S),\xi^2(S),\dots)$ is
defined as follows.  Set $\xi^k(S)=0$ if $k>\gamma_1$.
Next, let $k=\gamma_1$.  Define $\xi^k(S)$ as before, ignoring the
restriction on content.  Recall that for $j$ large, $V_k^j(S)$ stabilizes.
Let $\Sh$ be the tableau obtained by restricting 
this stable tableau to the numbers strictly greater than $m_k$,
where $m_k$ is the number of occurrences of the part $k=\gamma_1$
in $\gamma$.

Inductively define 
\begin{equation*}
  \xi^j(S)=\xi^j(\Sh) \qquad \text{ for $j<\gamma_1$}
\end{equation*}
for the tableau $\Sh$ of content 
$(0^{m_k},\gamma_{m_k+1},\gamma_{m_k+2},\dotsc,\gamma_n)$.

\begin{ex} Using the tableau $S$ in the previous example,
$\cmtype(S)$ is computed below.
\begin{equation*}
\begin{split}
  S &= \begin{matrix} 1&1&1&3&4&5\\ 2&2&2&4&5&6\\ 3&3&6& & & \\ 4&7&7& & &
  \end{matrix} \\
V_3(S) &= \begin{matrix} 1&1&1&6&7\\ 2&2&2&7& \\ 3&3&3& & \\
	4&4&4& & \\ 5&5& & & \\ 6& & & & \\ \end{matrix} \qquad
V^2_3(S) = \begin{matrix} 1&1&1\\ 2&2&2 \\ 3&3&3 \\
	4&4&4 \\ 5&5&  \\ 6&6& \\ 7&7&  \end{matrix}.
\end{split}
\end{equation*}
\begin{equation*}
  \Sh = \begin{matrix} 5&5\\ 6&6\\ 7&7 \end{matrix}.
\end{equation*}
We have $\cmtype(S)=((),(3),(2,2),(),\dots)$.  In this case
$\cmtype(S)=\xi(R)$.
\end{ex}

\begin{prop} \label{cat and type}
Let $R$ be such that $R_j\supseteq R_{j+1}$ for all $j$
and $S\in\CST(\gamma(R))$.  Then $\cmtype(S)\dom \xi(R)$
if and only if $S$ is $R$-column-catabolizable.
\end{prop}
\begin{proof} Let $\gamma=\gamma(R)$.
By definition $|\xi^i(S)|=|\xi^i(R)|$ for all $i$.
Let $k$ be maximal such that $\xi^k(S)$ is nonempty.
Write $M:=\xi^k(S)_1$ and $m:=\xi^k(R)_1$ so that
$R_1=(k^m)$.  Write $R'=(R_2,R_3,\dots)$ and
$\gamma'=(0^m,\gamma_{m+1},\dots,\gamma_n)$.

Note that $M\ge m$ since $\cmtype(S)\dom\xi(R)$
and $M$ and $m$ are the first parts of $\xi^k(S)$ and 
$\xi^k(R)$ respectively.  On the other hand,
if $S$ is $R$-column-catabolizable then $S|_{[m]}=Y_1$ and hence
$M\ge m$.  So without loss of generality it may be assumed that $M\ge m$
and $S|_{[m]}=Y_1$.

Divide $S$ into $Y_1$, $S_e$ and $S_w$ as in the definition
of $R$-column-catabolizability.
Then $S'=\colcat_{R_1}(S)=P(S_e S_w)$ and $V_k(S)=P(S_e S_w Y_1)$.

Because of the assumptions that $M\ge m$ and $S|_{[m]}=Y_1$,
it is enough to show that the following are equivalent:
\begin{enumerate}
\item $\cmtype(S) \dom \xi(R)$.
\item $\cmtype(S') \dom \xi(R')$.
\item $S'$ is $R'$-column-catabolizable.
\item $S$ is $R$-column-catabolizable.
\end{enumerate}
The equivalences $2\Leftrightarrow 3$ and 
$3\Leftrightarrow 4$ hold by induction and definition respectively.

$1 \Leftrightarrow 2$:  Note that
\begin{equation*}
  P(V_k^p(S') Y_1)=V_k^p(P(S' Y_1))=V_k^{p+1}(S)
\end{equation*}
for all $p\ge0$ since row insertion of $Y_1$ and $V_k$ commute on
on all column-strict tableaux of content $\gamma'$.
It follows from the definitions that
$\xi^j(S')=\xi^j(S)$ for $j<k$,
that the first part of $\xi^k(S')$ is $M-m+\xi^k(S)_2$,
and that the $j$-th part is $\xi^k(S)_{j+1}$ for all $j>1$.
Since $M\ge m$ it follows that $1 \Leftrightarrow 2$.
\end{proof}

\subsection{Image characterization revisited}

Proposition \ref{cat and type} and Theorem \ref{embed image} show that if
$R_j\supseteq R_{j+1}$ for all $j$, then $S\in \Img \theta_R$ if and only
if $\cmtype(S)\dom\xi(R)$.

For general $R$, in the notation of Proposition \ref{image},
one may reduce the test for $S\in\Img \theta_R$ to
that of $(\sigma(S))|_I\in \Img \theta_{R^k}$.  But the sequence of
rectangles $R^k$ has a rearrangement that is linearly ordered by
containment, since all of the rectangles in $R^k$ have $k$ columns.
Thus the catabolism multi-type can be used for each $k$.

\section{Transposition of poset structure}

The transposition maps are necessary to describe the relationship
between the embeddings $\theta$ and the cyclage-preserving
embeddings $\CST(\mu)\rightarrow \CST((1^n))$ given in \cite{La}.

\subsection{$R$-cyclage}
\label{cyclage sec}

The set $\LRT(R)$ has another structure as a graded poset
which shall be called the $R$-cyclage.  Its grading function, called
$\cocharge_R$, is complementary to $\charge_R$.
But it is not the dual of the $R$-cocyclage poset, but
rather the ``transpose'' of the $R^t$-cocyclage poset, where
$R^t$ is the sequence of rectangles given by the transposes
of the rectangles in $R$.

Let $b=\max_i \eta_i$ and $w=xu\in W(R)$ with $x$ a letter.  Denote
by $P(ux) <^R P(xu)$ for a covering relation for the $R$-cyclage poset;
this relation holds if the cell given by the difference
of the shapes of the tableaux $P(xu)$ and $P(u)$ is in a row
strictly south of the $b$-th.

Let us define a function $\cocharge_R:W(R)\rightarrow \NN$ in the
same way as $\charge_R$ except the function
$d_{R_1,R_2}:W((R_1,R_2))\rightarrow\NN$
is replaced by $\cod_{R_1,R_2}:W((R_1,R_2))\rightarrow\NN$,
where $\cod_{R_1,R_2}(w)$ is the number of cells of the shape of
$P(w)$ that are strictly below the $\max(\eta_1,\eta_2)$-th row.
Note that $d_{R_1,R_2}(w)+\cod_{R_1,R_2}(w) = | R_1 \cap R_2 |$.

\begin{thm} \label{cyclage theorem}
\begin{enumerate}
\item $\LRT(R)$ is a graded poset under the $R$-cyclage relation
with grading function $\cocharge_R$.
\item A tableau in $\LRT(R)$ is $R$-cyclage minimal if and only if
it has exactly $b$ rows.
\item $\tau_p$ is an isomorphism from $\LRT(R)$ under $R$-cyclage
to $\LRT(\tau_p R)$ under $(\tau_p R)$-cyclage.
\end{enumerate}
\end{thm}

The proof is similar to those of \cite[Theorems 19, 20]{S}.
The statistics $\charge_R$ and $\cocharge_R$ are complementary
in the following sense.
Let $r_{i,j}(R)$ be the number of indices $k$ such that
the rectangle $R_k$ contains the cell $(i,j)$.  Write
\begin{equation*}
  n(R) = \sum_{(i,j)} \binom{r_{i,j}(R)}{2}.
\end{equation*}

In the Kostka case where $\eta_i=1$ for all $i$,
$n(R)=n(\mu)=\sum_i (i-1)\mu_i$ as in \cite[I.1.5]{Mac},
and the maps $\charge_R$ and $\cocharge_R$ are the
charge and cocharge statistics of Lascoux and Sch\"utzenberger.

\begin{prop} \label{charge comp} For every word $w\in W(R)$,
$\charge_R(w) + \cocharge_R(w) = n(R)$.
\end{prop}
\begin{proof} Let $R=(R_1,\dots,R_t)$.  Let us write
$\chi(P)$ to give the value 1 if $P$ is true and $0$ otherwise.
Starting with the definition of $\charge_R$ in \cite[Section 3.2]{S},
we have
\begin{equation*}
\begin{split}
  &\,\,\charge_R(w) + \cocharge_R(w) = \dfrac{1}{t!}
  	\sum_{\sigma\in S_t} \sum_{i=1}^{t-1} (t-i)
		(d_{i,\sigma R}(\sigma w)+
		\cod_{i,\sigma R}(\sigma w)) \\
  &= \dfrac{1}{t!}
    	\sum_{\sigma\in S_t} \sum_{i=1}^{t-1} (t-i)
		| \sigma(R)_i \cap \sigma(R)_{i+1} |  \\
  &= \dfrac{1}{t!} \sum_{1\le k<l\le t} |R_k \cap R_l|
  	\sum_{i=1}^{t-1} (t-i) \sum_{\sigma\in S_t}
	\delta_{k,\sigma(i)}\delta_{l,\sigma(i+1)}  \\
  &= \dfrac{1}{t!} \sum_{1\le k<l\le t} |R_k \cap R_l|
  	\sum_{i=1}^{t-1} (t-i) (t-2)! 
   = \dfrac{1}{2} \sum_{1\le k<l\le t} |R_k \cap R_l| \\
  &= \dfrac{1}{2} \sum_{1\le k<l\le t} \sum_{(i,j)}
  	\chi((i,j)\in R_k) \chi((i,j)\in R_l) 
   = \dfrac{1}{2} \sum_{(i,j)} r_{i,j}(R) (r_{i,j}(R)-1) \\
  &= n(R).
\end{split}
\end{equation*}
\end{proof}

\subsection{Transpose bijection}

Let $R$ be a sequence of rectangles.  Recall the definition of the
subintervals $A_i$ from \cite[Section 2.2]{S}.
Let $R^t$ denote the sequence of transposes of the rectangles of $R$
and $A^t_i$ the subintervals for $R^t$.

Define the map $\tr_R:W(R)\rightarrow W(R^t)$ by
$w\mapsto w'$ where $w'$ is obtained from $w$ by
replacing the $c$-th copy (from the left) of the $r$-th largest
letter of $A_i$ by the $c$-th largest letter in $A^t_i$ for all $i$, $c$,
and $r$, and then reversing the resulting word.  To see that
$\tr_R(w)\in W(R^t)$, by restriction to the interval $A^t_i$ it is enough
to check it for the case $R=(R_1)$, in which case it follows from
the fact that $w\in W((R_1))$ if and only if $w$ is lattice of content
$R_1$.

Using the same words as above, one has a bijection
$\tr_R:\LRT(R)\rightarrow\LRT(R^t)$ that
transposes the shapes of tableaux.  This bijection on LR tableaux
appears in \cite{KS}.

\begin{ex} Let $R$ and $\la$ be as in the previous example.\newline
We have $R^t=((222),(222),(33))$, $A^t_1=[1,3]$,
$A^t_2=[4,6]$ and $A^t_3=[7,8]$, with tableaux
\begin{equation*}
  Y^t_1 = \begin{matrix} 1&1\\2&2\\3&3\end{matrix} \qquad
  Y^t_2 = \begin{matrix} 4&4\\5&5\\6&6\end{matrix} \qquad
  Y^t_3 = \begin{matrix} 7&7&7\\8&8&8\end{matrix}.
\end{equation*}
A tableau $T\in\LRT(\la;R)$ and $\tr_R(T)\in\LRT(\la^t;R^t)$ are given
below.
\begin{equation*}
T=\begin{matrix}
  1&1&1&3&3&5\\
  2&2&2&4&5&6\\
  3&4&6& & & \\
  4&7&7& & &
\end{matrix} \qquad
T^t =
  \begin{matrix}
    1&1&4&4\\2&2&5&7\\3&3&7&8\\5&6& & \\6&7& & \\8&8& &
  \end{matrix}
\end{equation*}
\end{ex}

\begin{thm} \label{poset transpose}
The map $\tr_R:\LRT(R)\rightarrow \LRT(R^t)$ is a bijection that transposes
the shapes of tableaux and gives a poset isomorphism from $R$-cocyclage to
$R^t$-cyclage that sends $\charge_R$ to $\cocharge_{R^t}$.
\end{thm}

\begin{rem} Suppose in Theorem \ref{poset transpose}
the $R$-cocyclage is replaced by the strong $R$-cocyclage
and $R^t$-cyclage by an obvious definition of strong $R^t$-cyclage.
The resulting theorem follows immediately from Propositions
\ref{chi trans} and \ref{switch trans}.
\end{rem}

Before proving Theorem \ref{poset transpose}, note that
it and Proposition \ref{charge comp} solve the
following conjecture of \cite{KS}.

\begin{thm} \label{poly transpose} For $R$ dominant and
$R'$ a dominant rearrangement of $R^t$,
\begin{equation*}
\K_{\la^t;R'}(q) = q^{n(R)} \K_{\la;R}(q^{-1})
\end{equation*}
\end{thm}

\begin{rem} \label{cyc std}
In the Kostka case (that is, $R_i=(\mu_i)$ for all $i$),
we have $n(R)=n(\mu):=\sum_i (i-1)\mu_i$, and Theorem \ref{poly transpose}
yields
\begin{equation*}
  \K_{\la^t;((1^{\mu_1}),(1^{\mu_2}),\dots)}(q) =
	q^{n(\mu)} \K_{\la,\mu}(q^{-1}) =: \tK_{\la,\mu}(q)
\end{equation*}
where the left hand side is the Poincar\'e polynomial and the
right hand side is the normalized Kostka-Foulkes polynomial.
The cyclage standardization map in \cite{La} is a cocharge-preserving
map from column-strict tableaux of content $\mu$ (where $\mu$ is a
partition of $n$) to standard tableaux (those of content $(1^n)$),
is precisely the map $\tr_{((1)^n)} \circ \theta_{R^t} \circ \tr_R$.
Note that $\tr_{((1)^n)}$ is just the
obvious transposition of standard tableaux
that have $n$ letters.  If $T$ is the input tableau and $S$ the
output tableau, then one has $\charge(S) = n((1^n))-(n(\mu)-\charge(T))$
by the properties of $\tr$ and $\theta$.  But this means
precisely that $\cocharge(T)=\cocharge(S)$.
\end{rem}

\begin{prop} \label{trans lemma} Let $w\in W(R)$.
\begin{enumerate}
\item[(T1)] $\tr_{R^t} \circ \tr_R$ is the identity on $W(R)$ and $\LRT(R)$.
\item[(T2)] $w$ is the row-reading word of a column-strict tableau
of (possibly skew) shape $D$ if and only if $\tr_R(w)$ is the
column-reading word of a column-strict tableau of shape $D^t$.
\item[(T3)] $\tr_R$ preserves Knuth equivalence.
\item[(T4)] $P(\tr_R(w))=\tr_R(P(w))$.
\item[(T5)] $Q(\tr_R(w))=\ev(Q(w))^t$ where $\ev$ means evacuation.
\end{enumerate}
\end{prop}
\begin{proof} (T1) follows easily from the definitions.
Let $w\in W(R)$ so that $\tr_R(w)\in W(R^t)$.
By Remark \ref{tr comp}, $\rev(\cstd(w))=\std(\tr_R(w))$ holds and this
equality characterizes $\tr_R(w)$.

By Lemma \ref{std lem} (see the appendix)
the maps $\cstd$ and $\std$ preserve the property of
being the row-reading word of a column-strict tableau of shape $D$,
send Knuth equivalence to Knuth equivalence, and 
preserve $Q$ tableaux.

Let $v$ be a standard word.  It is obvious that
$v$ is the row-reading word of a standard tableau of shape $D$,
if and only if $\rev(v)$ is the column-reading word of a standard
tableau of shape $D^t$.  This proves (T2).
It is well-known that $Q(\rev(v))=\ev(Q(v)^t$, proving (T5).
Reversal of standard words preserves Knuth equivalence, proving (T3)
and (T4).
\end{proof}

\begin{prop} \label{chi trans} For all $w\in W(R)$,
\begin{equation} \label{chi and trans}
  \tr_R(\chi_R(w)) = \chi_{R^t}^{-1}(\tr_R(w)).
\end{equation}
\end{prop}
\begin{proof} Let $A^t_j$ be the subalphabet for the
transposed rectangle $R^t_j$.  Now $\chi_R$ cyclicly rotates the
positions of the letters of $A_j$ in $w$ and $\tr_R$
sends the letters of $A_j$ to letters of $A_j^t$ occurring
in positions that are complementary with respect to the
length of the word $w$.  This given, it is clear that
the letters of $A^t_j$ occur in the same set of positions
in both words in \eqref{chi and trans}.  Therefore it suffices
to show that the restrictions of each of these words to $A^t_j$
coincide for all $j$.  Suppose the last letter of $w$ is in $A_j$.
Then the last letter of $\tr_R(w)$ is in
$A_j^t$.  By the definition of $\chi_R$ and $\tr_R$, the
restrictions of the two words to $A^t_i$ agree for all $i\not=j$.

Thus it may be assumed that $R=(R_1)$.  The $P$-tableaux of both
sides are equal since there is only one $R$-LR tableau for $R=(R_1)$.
For the standard row-insertion recording tableaux, recall from \cite{S}
the promotion operator $\pr_1$, evacuation $\ev$, and 
shape transposition $()^t$.  Then
\begin{equation*}
\begin{split}
  Q(\tr_R(\chi_R(w))) &= \ev(Q(\chi_R(w)))^t \\
  &= \ev(\pr_1(Q(w)))^t \\
  &= \pr_1^{-1}(\ev(Q(w))^t) \\
  &= Q(\chi_{R^t}^{-1}(\tr_R(w))).
\end{split}
\end{equation*}
These equalities of rectangular standard tableaux hold by
(T5), \cite[Proof of Proposition 15]{S}, the fact that
$\ev \pr_1 = \pr_1^{-1} \ev$ on standard rectangular tableaux,
and the commutation of $\ev$ and $\pr_1$ with transpose.
\end{proof}

The transpose map is also compatible with rectangle-switching.

\begin{prop} \label{switch trans} For all $w\in W(R)$,
\begin{equation} \label{switch and trans}
  \tr_{\tau_p R}(\tau_p(w)) = \tau_p(\tr_R(w)).
\end{equation}
\end{prop}
\begin{prop} Let $B=A_p\cup A_{p+1}$.  From the definitions,
it is easy to check that at positions in the two words in
\eqref{switch and trans} that contain a letter not in $B$, the two
words coincide.  Thus it may be assumed that $R=(R_1,R_2)$ and $p=1$.
The $P$ tableaux of the words are equal,
\begin{equation*}
  P(\tr_{\tau_p R}(\tau_p(w)))=\tr_{\tau_p R}(\tau_p(P(w)))=
  \tau_p(\tr_R(P(w))) = P(\tau_p(\tr_R(w))),
\end{equation*}
by (T4), \cite[Theorem 9 (A3)]{S}, and the fact that
$V_{R_1^t} \otimes V_{R_2^t}$ is multiplicity-free.
The $Q$ tableaux also coincide:
\begin{equation*}
Q(\tr_{\tau_p R}(\tau_p(w)))=\ev(Q(\tau_p(w)))^t =
  \ev(Q(w))^t = Q(\tr_R(w)) = Q(\tau_p(\tr_R(w))),
\end{equation*}
by (T5) and \cite[Theorem 9 (A4)]{S}.
\end{prop}

Proof of Theorem \ref{poset transpose}:
\begin{proof} 
Let $S<_R T$ be an $R$-cocyclage covering relation in $\LRT(R)$.
Let $s$ be the corner cell of the shape $\la$ of $T$, $x$ be the letter
and $U$ the column-strict tableau of shape $\la-s$, such that
$T=P(Ux)$ and $S=P(\chi_R(Ux))=P((w_0^R x)(w_0^R U))$.

Write $S'=\tr_R(S)$, $T'=\tr_R(T)$, and $x'u=\tr_R(Ux)$ where $x'$ is a
letter.  It must be shown that $S' <^{R^t} T'$.  We have
\begin{equation*}
\begin{split}
  T'&=\tr_R(T)=\tr_R(P(Ux))=P(\tr_R(Ux))=P(x'u) \\
  S'&=\tr_R(S)=\tr_R(P(\chi_R(Ux)))\\
  &=P(\tr_R(\chi_R(Ux)))=P(\chi_{R^t}^{-1}(\tr_R(Ux)))
\end{split}
\end{equation*}
by (T4) and Proposition \ref{chi trans}.  By the definition
of $R$-cocyclage, it is enough to show
that $P(u)$ and $U$ have transpose shape.  Now
$Ux$ fits the skew shape $D \otimes (1)$ where $D$ is the
shape of the tableau $U$.  By (T2) $\tr_R(Ux)=x'U$ fits
the skew shape $(D \otimes (1))^t=(1)\otimes D^t$.
In fact by definition $u$ is the column-reading word
of a column-strict tableau $U'$ of shape $D^t$.
Thus $P(u)=U'$ and $U$ have transpose shape,
and $\tr_R$ sends the $R$-cocyclage to $R^t$-cyclage.

To show that the transpose map sends $\charge_R$ to $\cocharge_{R^t}$, 
by the above poset isomorphism it is enough 
to show that this holds for $R$-cocyclage minimal
tableaux.  In light of Proposition \ref{switch trans}
it may be assumed that $R$ is dominant.

Let $Y_1$ and $Y_1^t$ denote the first Yamanouchi tableaux
for $R$ and $R^t$ respectively, $\Rhat=(R_2,R_3,\dots)$ and
$\Rhat^t=(R_2^t,R_3^t,\dots)$.  Let $T\in\LRT(\la;R)$ be $R$-cocyclage
minimal and let $T$ consist of $Y_1$ sitting atop $\Th\in\LRT(\Rhat)$.
Since $\tr_R$ is shape-transposing it follows that
$\tr_R(T)$ is $R^t$-cyclage minimal and therefore consists
of $Y_1^t$ sitting to the left of a tableau $\Th'\in\LRT(\Rhat^t)$.
By the definition of $\tr_R$ it is easy to see that
$\Th'=\tr_{\Rhat}(\Th)$.  Then by induction
\begin{equation*}
\begin{split}
  \charge_R(T)&=\charge_{\Rhat}(\Th)=
  \cocharge_{\Rhat^t}(\tr_{\Rhat}(\Th))\\
  &= \cocharge_{\Rhat^t}(\tr_R(T)-Y_1^t)=
  \cocharge_{R^t}(\tr_R(T)).
\end{split}
\end{equation*}
\end{proof}

\section{Appendix}

The purpose of this section is to describe a generalization of
Schensted's standardization map 
that has enough flexibility for use with the bijection $\tr_R$
on $R$-LR words and tableaux.  First we recall the original
standardization map of Schensted \cite{Sch}.

\begin{prop} \label{std image}
Let $\std$ be the map (Schensted's standardization)
from words of length $N$ to words of content $(1^N)$
such that $\std(w)_i < \std_(w)_j$ if $w_i<w_j$, or $w_i=w_j$ and $i<j$.
The the following are equivalent for $\alpha=(\alpha_1,\alpha_2,\dots)$.
\begin{enumerate}
\item There is a (necessarily unique) word $w$ of content $\alpha$ such that
$v=\std(w)$.
\item For every value $i$ such that $i+1$ precedes $i$ in $v$,
$i$ must have the form $i=\alpha_1+\dots+\alpha_r$ for some $r$.
\end{enumerate}
\end{prop}

Let $A=A_1+A_2+\dots+A_t$ be a segmentation of an alphabet $A$,
that is, a set partition of $A$ into subintervals $A_i$ such that
if $i<j$ then $x<y$ for all $x\in A_i$ and $y\in A_j$.

Let $R=(R_1,R_2,\dots,R_t)$ be a sequence of partitions
and $Y=(Y_1,\dots,Y_t)$ a sequence of column-strict tableaux
such that $Y_i$ has shape $R_i$ and letters in $A_i$ for all $i$.
Denote by $W(R,Y)$ the set of words $w$ in the alphabet $A$
such that $P(w|_{A_i})=Y_i$ for all $i$.  The set of $R$-LR words
is a special case of this construction.

Let $B=B_1+B_2+\dots+B_t$ be a segmentation of the alphabet $B$
and $Z=(Z_1,\dots,Z_t)$ such that $Z_i$ is column-strict of
shape $R_i$ in the alphabet $B_i$ for all $i$.  Due to the bijectivity
of the RS correspondence, there is a unique bijection
$\std_Y^Z:W(R,Y)\rightarrow W(R,Z)$ such that:
\begin{enumerate}
\item For all $i$, the positions of the letters of $B_i$ in $\std_Y^Z(w)$
are equal to those of $A_i$ in $w$.
\item $P(\std_Y(w)|_{B_i})=Z_i$ for all $i$.
\item $Q(\std_Y(w)|_{B_i})=Q(w|_{A_i})$ for all $i$.
\end{enumerate}

For a partition $\la$ of $N$, define the rowwise tableau of shape
$\la$ to be the standard tableau of shape $\la$ in which the
first row is comprised of the first $\la_1$ letters in the
interval $[N]$, the second row is comprised of the next $\la_2$
letters in $[N]$, etc.  Define the columnwise tableau of shape
$\la$ in the same way except that rows are replaced by columns and
$\la_i$ by $\la^t_i$.

Let $R=(R_1,\dots,R_t)$ be a sequence of partitions,
$B=[N]$ where $N=\sum_i |R_i|$, and $B=B_1+B_2+\dots+B_t$ the segmentation
of $B$ such that $|B_i|=|R_i|$.  For all $i$, let $Z_i$ be the
rowwise standard tableau of shape $R_i$ in the alphabet $B_i$.
Then the map $\std_Y^Z$ coincides with Schensted's standardization map
$\std$ (see Lemma \ref{std lem}).  If $Z_i$ is taken to be
the columnwise standard tableau of shape $R_i$ in the alphabet $B_i$
then write $\cstd$ for the map $\std_Y^Z$.

\begin{rem} \label{tr comp} Now let us require that
each of the $R_i$ be rectangles.  Let $w\in W(R)$,
so that $\tr_R(w)\in W(R^t)$.
Then by construction, $\rev(\cstd(w))=\std(\tr_R(w))$.  Moreover,
this condition determines $\tr_R(w)$ uniquely in light of
Proposition \ref{std image}.  It should be mentioned that the map
$\tr_R$ can be defined for arbitrary sequences of partitions,
or more generally, for skew shapes,
using Zelevinsky's definition of a picture \cite{Z}.
Even more generally, in \cite{RS} there is a transpose
construction that generalizes all of the above.

In the present case it is possible to use a direct relabeling,
as opposed to the Schensted constructions that must be used in \cite{RS}.
\end{rem}

Here are the main properties of the generalized standardization maps.
A nontrivial construction (a crystal lowering operator)
is used in the proof.

\begin{lem} \label{std lem} In the above notation,
\begin{enumerate}
\item $w\in W(R,Y)$ is the row-reading word of a column-strict
tableau of the (skew) shape $D$ if and only if
$\std_Y^Z(w)\in W(R,Z)$ is.  In this case, if
$w=\word(T)$ where $T$ is column-strict of shape $D$, then
write $\std_Y^Z(T)$ for the column-strict tableau of shape $D$
whose word is $\std_Y^Z(w)$.
\item If $v\Knuth w$ where $v,w\in W(R,Y)$ then
$\std_Y^Z(v)\Knuth \std_Y^Z(w)$ in $W(R,Z)$.
\item $P(\std_Y^Z(w))=\std_Y^Z P(w)$.
\item $Q(\std_Y^Z(w))=Q(w)$.
\end{enumerate}
\end{lem}
\begin{proof} Suppose one of the tableaux $Y_i$ is not the
unique column-strict tableau of shape and content $R_i$
in the alphabet $A_i$.  Then $Y_i$ admits a crystal lowering operator
$e_r$, say, where $r\in A_i$.  Then $e_r$ induces the bijection
$\std_Y^{Y'}$ from $W(R,Y)$ to $W(R,Y')$
where $Y'_j=Y_j$ for $j\not=i$ and $Y'_i = e_r Y_i$.
But it is well-known that $e_r$ satisfies the above properties.
By applying this process to alter $Y$ and $Z$, it may be assumed
that for all $i$, $Y_i$ and $Z_i$ are the unique column-strict
tableaux of shape and content $R_i$ in the alphabets $A_i$ and $B_i$
respectively.  In this case $\std_Y^Z$ is a trivial relabeling and the
above properties are obviously satisfied.
\end{proof}

\end{document}